\documentclass[11pt,a4paper]{article}
\usepackage[T1]{fontenc}
\usepackage[utf8]{inputenc}
\usepackage[british]{babel}
\usepackage{lmodern}
\usepackage[a4paper,margin=1in]{geometry}
\usepackage{amsmath,amssymb,amsthm}
\usepackage{booktabs}
\usepackage{microtype}
\usepackage{xcolor}
\usepackage{xurl}
\definecolor{linkblue}{RGB}{0,70,140}
\usepackage[colorlinks=true,linkcolor=linkblue,citecolor=linkblue,urlcolor=linkblue]{hyperref}
\hypersetup{
  pdftitle={Big Ramsey degrees and structures in relational languages},
  pdfauthor={Yifan Hu},
  pdfkeywords={big Ramsey structures, infinite languages, Fra\"iss\'e structures, free amalgamation}
}
\newcommand{\DOI}[1]{\href{https://doi.org/#1}{\nolinkurl{doi:#1}}}
\numberwithin{equation}{section}
\newtheorem{theorem}{Theorem}[section]
\newtheorem{lemma}[theorem]{Lemma}
\newtheorem{proposition}[theorem]{Proposition}
\newtheorem{corollary}[theorem]{Corollary}
\theoremstyle{definition}
\newtheorem{definition}[theorem]{Definition}
\newtheorem{remark}[theorem]{Remark}
\DeclareMathOperator{\Emb}{Emb}
\DeclareMathOperator{\Age}{Age}
\DeclareMathOperator{\Aut}{Aut}
\DeclareMathOperator{\NF}{NF}
\DeclareMathOperator{\val}{val}
\DeclareMathOperator{\Str}{Str}

\newcommand{\N}{\omega}
\newcommand{\res}{\mathbin{\upharpoonright}}
\newcommand{\cat}{\mathbin{{}^\frown}}
\newcommand{\brd}{T_{\mathrm{emb}}}
\newcommand{\crd}{T_{\mathrm{copy}}}
\newcommand{\Sym}{\mathfrak S}

\title{Big Ramsey degrees and structures in relational languages}
\author{Yifan Hu\\[2pt]\normalsize Shanghai, China}
\date{}

\begin{document}
\maketitle

\begin{abstract}
We establish a precompact transfer theorem for big Ramsey structures and a uniform envelope theorem for valuation hypergraphs. These yield big Ramsey structures for generic injective relational structures in languages with finitely many symbols of each positive arity, answering the existence part of~\cite[Problem~6.3.2]{BinaryExact}, and refute~\cite[Conjecture~4.7]{InfiniteLanguages}. We also construct a finitely constrained free amalgamation class in a finite language of maximum arity three whose Fraïssé limit has infinite singleton big Ramsey degree, giving a negative answer to the existence question in~\cite[Problem~6.3.3]{BinaryExact}.
\end{abstract}


\section{Introduction}\label{sec:intro}

Ramsey's 1930 theorem asserts that a finite colouring of the $k$-element subsets
of an infinite set becomes constant on an infinite subset~\cite{Ramsey}.
For a countable relational structure $K$, the corresponding question is
whether the subset can be required to induce a copy of $K$. In general,
one colour no longer suffices. The big Ramsey degree of a finite structure
$A$ in $K$ measures how many colours may be unavoidable. It asks for a
bound depending on $A$ and $K$, but not on the colouring or the number of
colours. Determining these degrees therefore identifies precisely how
much of the infinite Ramsey theorem survives the requirement that the
ambient structure be preserved.

Devlin's 1979 calculation for the rationals~\cite{Devlin} and Milliken's
partition theorem for trees~\cite{Milliken} established two features of
the subject: finite tree configurations can describe unavoidable
colourings, and a tree theorem can show that these are the only
obstructions. For the Rado graph, Sauer proved finite upper bounds and
Laflamme, Sauer and Vuksanovic identified optimal canonical
partitions~\cite{Sauer,LSV}. Here a canonical partition is a fixed finite
partition of the copies of $A$ such that every finite colouring becomes
constant on its parts in a suitable copy of $K$, while every part occurs
in every copy. Thus the types used in an upper-bound proof have
structural meaning as well as a counting role.

Zucker's 2019 theory of big Ramsey structures asks that these descriptions
hold simultaneously for all finite $A$~\cite{ZuckerDynamics}. Such a
structure is an expansion of $K$ in which the induced finite expansions
record the optimal unavoidable distinctions. The expansion types must
agree under restriction to smaller structures. Finiteness of the degrees
alone does not provide this compatibility. Big Ramsey structures also
connect these partition problems with completion flows of automorphism
groups; see~\cite{ZuckerDynamics,HZ}. In this paper the relevant point is
combinatorial: one expansion should account for all finite colourings,
rather than supply unrelated optimal partitions for individual $A$.

The development beyond unrestricted binary structures brought the
effect of forbidden configurations into focus. Dobrinen obtained finite
big Ramsey degrees for the Henson graphs~\cite{DobrinenTriangleFree,DobrinenHenson}.
Zucker extended finiteness to finitely constrained binary free
amalgamation classes~\cite{ZuckerBinary}, and Balko, Chodounsk\'y,
Dobrinen, Hubi\v cka, Kone\v cn\'y, Vena and Zucker subsequently
characterized their exact degrees and constructed big Ramsey
structures~\cite{BinaryExact}.

The generic $3$-uniform hypergraph provides the first test case beyond
binary languages in that programme. Balko, Chodounsk\'y, Hubi\v cka,
Kone\v cn\'y and Vena proved that it has finite big Ramsey degrees by
coordinating two valuation trees~\cite{TernaryFinite}. Problem~6.3.1
of~\cite[p.~2144]{BinaryExact} asks for the next step: a characterization
of the exact degrees and a big Ramsey structure realizing them
simultaneously. Hubi\v cka, Kone\v cn\'y and Zucker have since announced
a characterization by hypergraph diaries, with Ramsey and unavoidability
theorems that also yield a big Ramsey structure
\cite[Theorems~1.4, 3.2 and~3.3]{TernaryDiaries}.

Balko et al. formulate two natural generalizations
\cite[Section~6.3]{BinaryExact}. Problem~6.3.2 asks for
the corresponding degree characterization and big Ramsey structure for
unrestricted relational structures in languages with finitely many
symbols of each arity; both the language and its arities may be unbounded.
Problem~6.3.3 asks the same questions when the language is finite and
finitely many irreducible structures are forbidden. These two directions
organize the present paper: we establish the existence assertion in the
first setting and give a negative instance of that assertion in the second.

\subsection{Unrestricted structures}

For the unrestricted direction, Braunfeld, Chodounsk\'y, de Rancourt,
Hubi\v cka, Kawach and Kone\v cn\'y extended the valuation construction
to languages with finitely many relations of each arity, proving finite
big Ramsey degrees throughout this range~\cite{InfiniteLanguages}.
Our existence result concerns the generic limit of all finite injective
$L$-structures in each such language.

Except in Theorem~\ref{thm:transfer} and its proof, languages have positive
arities and relations hold only on tuples of distinct entries. We call
this the \emph{injective} convention;
symmetry is not assumed. As in~\cite{InfiniteLanguages}, the convention
is essential for unbounded arity: repetitions can encode infinitely many
independent relations on a fixed finite set and prevent a countable
universal structure. With finitely many symbols per arity, the injective
convention leaves only finitely many structures of each finite size.

For the existence argument we use a countable extension of the upper-bound
transfer and finite-palette recurrence arguments of
Hubi\v cka--Zucker~\cite[Theorem~3.10, Fact~3.4 and the proof of
Theorem~3.18]{HZ}. We record the form needed when the initial pullback is
precompact and give its compatibility argument in
Section~\ref{sec:extraction}. For an expansion $X$ of $K$ and a finite
$A\in\Age(K)$, write
$\mathcal P_A(X)=\{X\cdot a:a\in\Emb(A,K)\}$, where $X\cdot a$ is the
pullback to the fixed underlying set of $A$. The expansion is precompact
if every $\mathcal P_A(X)$ is finite.

\begin{theorem}[Transfer through a precompact expansion]\label{thm:transfer}
Let $K$ and $M$ be countably infinite structures in the same relational
language, and let $N$ be a relational expansion of $M$. Suppose there
are embeddings $\phi:K\to M$ and $\theta:M\to K$ such that
$X=N\cdot\phi$ is precompact and $\brd(B,N)=1$ for every finite
$B\in\Age(X)$. Then $K$ admits a big Ramsey structure. For every finite
$A\in\Age(K)$,
\begin{equation}\label{eq:transfer-min}
\brd(A,K)=\min_{\eta\in\Emb(K)}|\mathcal P_A(X\cdot\eta)|.
\end{equation}
If $\Age(K)$ has only finitely many isomorphism types of each finite size,
the expansion can be obtained by adding finitely many injective relation
symbols of each arity.
\end{theorem}

No homogeneity assumption is required of $K$ or $M$. The theorem and its
proof allow arbitrary arities and repeated entries in the original
relations; only the added relations are required to be injective.
Applying it to the valuation structures constructed below gives the
following result.

\begin{corollary}\label{thm:brs}
Let $L$ be a relational language with finitely many symbols of each
positive arity, and let $K$ be the Fra\"iss\'e limit of all finite
injective $L$-structures. Then $K$ admits a big Ramsey structure obtained
by adding finitely many injective relation symbols of each arity.
\end{corollary}

The corollary answers the existence part of Problem~6.3.2 under this
convention. The same existence question is explicitly restated as
Problem~8.19 of Hubi\v cka--Kone\v cn\'y~\cite{HKClassification}.
Proposition~\ref{prop:degrees} also gives an implicit minimum
formula using a fixed normal-form colouring. Determining which collections
of normal forms realize the minimum is a further question.
The diary announcement also discusses extensions to finite relational
languages~\cite[Section~4]{TernaryDiaries}. Corollary~\ref{thm:brs}
includes unbounded arity; its minimum formula is an existence description,
rather than an effective diary enumeration.

The combinatorial ingredient is a uniform envelope theorem for the
valuation hypergraphs of~\cite[Definition~3.1]{InfiniteLanguages}. Suppose that
$L$ has $n_r<\infty$ relation symbols of arity $r\geq2$. An
$L$-hypergraph has symmetric injective relations, and each finite set
supports at most one relation. Put $\sigma_j=n_{j+1}+2$, and let
$G_\sigma$ be the valuation hypergraph recalled in
Section~\ref{sec:envelopes}.

\begin{theorem}\label{thm:envelopes}
Every $L$-hypergraph $H$ enumerated by $\N$ admits an embedding
$\phi:H\to G_\sigma$ preserving the enumeration such that, for every
$k\geq1$ and every $S\in[H]^k$, the set $\phi[S]$ has a valuation
envelope of height at most
\begin{equation}\label{eq:R}
R(k)=2k+(k4^k+1)^2.
\end{equation}
The same embedding works for all $k$, and the bound is independent of
$\sigma$ and $H$.
\end{theorem}

The quantifier order matters. The earlier upper-bound argument chooses
an embedding after the size $k$ has been fixed. When the numbers $n_r$
are unbounded, Conjecture~4.7 of~\cite{InfiniteLanguages} asserts that
every single embedding has some fixed finite size for which envelope
heights are unbounded. Theorem~\ref{thm:envelopes} uses the same ambient
$G_\sigma$ and the same notion of envelope, and gives the opposite
quantifier pattern.

\begin{corollary}\label{cor:counterexample}
Conjecture~4.7 of~\cite{InfiniteLanguages} is false.
\end{corollary}

The construction allows slices of order $m$ to meet at an original
vertex $n$ when $n<m$. For a $k$-element set $S$, the original vertices
subsequently added by closure therefore lie in
\[
S\cup\{0,\ldots,k-1\}.
\]
This finite set controls all later closure operations. Coordinated tree
embeddings preserve the resulting finite types, and the
finite-dimensional form of Milliken's theorem supplies their Ramsey
property. To obtain one expansion, we construct finite-stage expansions
on the same underlying set and prove that each extends its predecessor.
Theorem~\ref{thm:transfer} isolates this compatibility argument from the
envelope bound.

\subsection{Finitely many forbidden structures}

For a finite relational language $L$ and a finite family $\mathcal F$ of
finite irreducible $L$-structures, let $\operatorname{Forb}_e(\mathcal F)$
be the class of finite injective $L$-structures omitting induced copies
of members of $\mathcal F$. Irreducibility means that every pair of
distinct vertices occurs together in some positive relation tuple.
This condition ensures free amalgamation. Problem~6.3.3
of~\cite[p.~2145]{BinaryExact} asks for the big Ramsey degrees of the
Fra\"iss\'e limits of these classes and asks whether they admit big
Ramsey structures. Our second result gives a negative instance of the
existence question.

\begin{theorem}\label{thm:negative}
There are a finite relational language $L_*$ of maximum arity three and
a finite family $\mathcal F_*$ of finite irreducible injective
$L_*$-structures, each with three or four vertices, such that
$\operatorname{Forb}_e(\mathcal F_*)$ has a countably infinite
Fra\"iss\'e limit $K_*$ and a one-point substructure $\delta$ satisfying
\[
\brd(\delta,K_*)=\infty.
\]
In particular, $K_*$ has no big Ramsey structure.
\end{theorem}

The language, forbidden family and singleton in this statement are
fixed. For every $q\geq2$ there is a $q$-colouring of the copies of
$\delta$ such that every copy of $K_*$ contains all $q$ colours. Thus
finite language and finite constraints, even with free amalgamation,
do not ensure finite big Ramsey degrees once ternary relations are
allowed. Theorem~\ref{thm:negative} does not classify the classes in
Problem~6.3.3, nor does it decide the corresponding question for a
specified class of hypergraphs.

The counterexample also gives a negative answer to the finite-language,
finitely-forbidden strong-amalgamation subquestion of Problem~6.1 in
Dobrinen's ICM report~\cite{DobrinenICM}. Indeed, the class constructed here
has free amalgamation, and hence strong amalgamation. It likewise answers
negatively the universal finiteness question in Question~8.4 of
Coulson--Dobrinen--Patel~\cite{CDPQuestion}.
The additional non-$3$-irreducibility hypothesis in that question is met:
the family $\mathcal F_*$ contains an $E$-complete forbidden triple with
all ternary relations false, obtained by violating
\eqref{neg:states}. Such a triple is irreducible but not $3$-irreducible.

The lower-bound method is Sauer's persistent-partition argument:
finite antichains of ages of one-point realization sets give persistent
vertex partitions~\cite[Lemma~5.1 and Theorem~5.1]{SauerCanonicalVertex2003};
see also~\cite[Sections~1, 9 and~10]{SauerHomogeneous2020}.
Our finite presentation realizes arbitrarily large antichains in one
language and one forbidden class. For each finite loopless directed graph
$H$, it provides a finite parameter set and a realized one-point type
with realization set $P_H$, such that
\[
\Age(K_*[P_H])\subseteq\Age(K_*[P_J])
\quad\Longleftrightarrow\quad
H\longrightarrow J,
\]
where the arrow denotes a directed graph homomorphism. Directed cycles
of distinct prime lengths give arbitrarily large antichains. Applying
Sauer's persistence argument then gives infinite singleton degree. We
include the short persistence proof for the distinguished singleton and
prove the finite presentation and its age characterization in full.

Together the results distinguish two sources of difficulty beyond the
binary theory. In the unrestricted case, unbounded arity requires
uniform control of finite closures and compatibility of optimal
expansions. Finite forbidden configurations can instead create an
obstruction already at the level of vertices. The first construction
provides the needed compatibility; the second exhibits such an
obstruction explicitly.
\section{Embeddings, expansions and trees}\label{sec:prelim}

All structures are countable and relational, and embeddings preserve and
reflect relations. Write $\Emb(A,B)$ for their embeddings and $\Emb(K)$ for
$\Emb(K,K)$. An enumeration will be used for coding, but embeddings in
Corollary~\ref{thm:brs} need not preserve it. For finite $A\in\Age(K)$, the
embedding big Ramsey degree $\brd(A,K)$ is the least positive integer $d$
such that every finite colouring of $\Emb(A,K)$ uses at most $d$ colours on
$g\circ\Emb(A,K)$ for some $g\in\Emb(K)$. If no finite $d$ works, put
$\brd(A,K)=\infty$. Write
\[
 \binom KA=\{B\subseteq K:K\res B\cong A\}.
\]
The copy big Ramsey degree $\crd(A,K)$ is defined in the same way, using
finite colourings of $\binom KA$ and their restrictions to $\binom{g[K]}A$;
again the value is $\infty$ when no finite bound exists. We use the standard
conversion between the two conventions~\cite[Section~2]{HZ}:
\begin{equation}\label{eq:degree-conventions}
\brd(A,K)=|\Aut(A)|\,\crd(A,K).
\end{equation}
For singletons the two degrees therefore agree.

If $X$ expands $K$ and $e:A\to K$ is an embedding, let $X\cdot e$ be its
pullback to the fixed underlying set of $A$. Put
\[
\mathcal P_A(X)=\{X\cdot e:e\in\Emb(A,K)\}.
\]
The expansion is \emph{precompact} if these sets are finite, and
\emph{unavoidable} if $\mathcal P_A(X\cdot g)=\mathcal P_A(X)$ for every
$A$ and $g\in\Emb(K)$. It is a \emph{big Ramsey structure} if it is
unavoidable and $|\mathcal P_A(X)|=\brd(A,K)<\infty$ for every $A$.
This is the embedding convention of~\cite[Definition~3.5]{HZ}, equivalent
to the usual formulation with expansion isomorphism types and automorphism
factors~\cite{ZuckerDynamics}.

We say that a countable structure $N$ has the \emph{infinite Ramsey property}
if every finite colouring of $\Emb(B,N)$, with $B\in\Age(N)$ finite, is
constant on $z\circ\Emb(B,N)$ for some $z\in\Emb(N)$. An expansion $F$ of
$K$ is \emph{recurrent} if every $g\in\Emb(K)$ admits $h\in\Emb(K)$ such
that $gh\in\Emb(F)$~\cite[Definition~3.3(4)]{HZ}. These definitions
will also be used for finite-stage expansions. The existence and extension
properties of Fra\"iss\'e limits are used in their standard form~\cite{Fraisse}.

A strong subtree of a rooted tree is a rooted, balanced, meet-closed subtree
whose levels lie on ambient levels and which, above every nonterminal node,
has exactly one next-level node in each immediate ambient successor cone.
A strong vector subtree has a common level set in all coordinates. Let
$\Str_h(T_0,\ldots,T_{d-1})$ denote the strong vector subtrees of height $h$.
Height counts levels of the subtree. We use the following standard form of
Milliken's theorem; see~\cite{Milliken} and~\cite[Theorem~2.2]{InfiniteLanguages}.

\begin{theorem}[Milliken]\label{thm:milliken}
Let $d,h\geq1$ be integers, and let $T_0,\ldots,T_{d-1}$ be rooted, finitely branching
trees with no terminal nodes. Every finite colouring of
$\Str_h(T_0,\ldots,T_{d-1})$ is constant on
$\Str_h(S_0,\ldots,S_{d-1})$ for some infinite strong vector subtree
$(S_i)_{i<d}$.
\end{theorem}

All applications below have finite dimension.
\section{Finite closures in valuation trees}\label{sec:source}

Fix a sequence $\sigma=(\sigma_j)_{j\geq1}$ of positive finite integers.
For $i,n<\omega$, let $T_i(n)$ consist of the height-tagged functions on
nonempty decreasing tuples below $n$ such that
\[
f(x_0,\ldots,x_{r-1})<\sigma_{i+r}.
\]
Set $T_i=\bigcup_nT_i(n)$, ordered by restriction, and write $|f|=n$ for
$f\in T_i(n)$. These are the valuation trees of~\cite[Section~3]{InfiniteLanguages}.
For $\ell<|f|$, its full slice is
\[
f_\ell\in T_{i+1}(\ell),\qquad
f_\ell(\bar y)=f(\ell,\bar y).
\]
The scalar $f(\ell)$ is separate. A successor of $f\in T_i(n)$ is uniquely
determined by this scalar at $n$ and its full $n$-slice in $T_{i+1}(n)$.
The trees have finite levels, since
\[
|T_i(n)|=\prod_{r=1}^n\sigma_{i+r}^{\binom nr},
\]
and have no terminal nodes.

\begin{definition}\label{def:closure}
For a nonempty finite $X\subseteq T_0$, let $C(X)=(C_i(X))_{i<\omega}$
be the least family containing $X$ in coordinate zero and closed under
same-coordinate meets, restrictions to all retained levels, and full slices
at all retained smaller levels. Here a level is retained if it is the height
of a node in any coordinate. Write $L(X)$ for the retained level set.
No zero nodes or other constants are adjoined.
\end{definition}

For $L(X)=\{\ell_0<\cdots<\ell_{h-1}\}$ define the compression of
$f\in C_i(X)$, $|f|=\ell_j$, by
\begin{equation}\label{eq:compression}
\widehat f(a_1,\ldots,a_r)=f(\ell_{a_1},\ldots,\ell_{a_r}),
\qquad j>a_1>\cdots>a_r.
\end{equation}

\begin{lemma}\label{lem:compression}
The closure $C(X)$ is finite. Compression identifies its labelled diagram
of nodes and ordered retained heights, including meets, restrictions,
full slices and scalar values at retained levels, with $C(\widehat X)$. Its level set becomes $\{0,\ldots,h-1\}$, and
\begin{equation}\label{eq:triangle}
f\in C_i(X),\quad |f|=\ell_j\quad\Longrightarrow\quad i+j<h.
\end{equation}
Compression is idempotent and commutes with passing to a subtuple, when the
subtuple is compressed again using its own closure.
\end{lemma}

\begin{proof}
Put $N=\max\{|f|:f\in X\}$. Induction on a derivation in $C(X)$ gives
\begin{equation}\label{eq:closure-invariants}
 f\in C_i(X)\ \Longrightarrow\
 i+|f|\leq N\quad\text{and}\quad
 |L(X)\cap(|f|,\infty)|\geq i.
\end{equation}
For the second assertion, a slice adds its parent's height to the higher
levels; a meet or restriction retains the higher levels of either parent.
The first assertion gives finiteness, and the second gives~\eqref{eq:triangle}.

Let $\rho(\ell_j)=j$. Directly from~\eqref{eq:compression},
\[
 \widehat{f\res\ell_j}=\widehat f\res j,\qquad
 \widehat{f_{\ell_j}}=(\widehat f)_j,\qquad
 \widehat f(j)=f(\ell_j)
\]
whenever the expressions are defined. We prove injectivity and preservation
of meets simultaneously by induction on the maximum ambient height of two
nodes. The comparable case follows from restrictions. Otherwise, put
$c=|f\wedge g|$. Since $c\in L(X)$, there are two cases:
\[
\begin{array}{ll}
 f(c)\ne g(c)
   &\Longrightarrow\widehat f(\rho(c))\ne\widehat g(\rho(c)),\\[2pt]
 f_c\ne g_c
   &\Longrightarrow\widehat{f_c}\ne\widehat{g_c}
    \Longrightarrow(\widehat f)_{\rho(c)}\ne(\widehat g)_{\rho(c)}.
\end{array}
\]
The second line uses induction, since $|f_c|=|g_c|=c$.
There is no earlier disagreement, and hence
\[
 \widehat f\wedge\widehat g=\widehat{f\wedge g}.
\]
Thus the compressed image is closed under all generating operations.
Minimality gives $C(\widehat X)\subseteq\widehat{C(X)}$; induction on
derivations gives the reverse inclusion. The same argument applied to
$Y\subseteq X$ proves compatibility with subtuples. Compression on
consecutive levels is the identity.
\end{proof}

For a finite nonempty label set $U$, the \emph{normal form} $\NF(\bar a)$ of an ordered
tuple of distinct vertices in $U\times T_0$ records its labels and the
compressed projections, in the original order. Distinct vertices may have
the same projection. We name every such normal form by an injective relation
symbol of the corresponding arity. Denote this countable expansion of
$U\times T_0$ by $N_{\sigma,U}$. Lemma~\ref{lem:compression} implies that
the normal form of a tuple determines that of each ordered subtuple.

\subsection{Coordinated embeddings}

Let $(S_i)_{i<\omega}$ be strong subtrees of $(T_i)_{i<\omega}$ with common
level sequence $e(0)<e(1)<\cdots$. There are maps $\Psi_i:T_i\to S_i$
defined by the simultaneous recursion of~\cite[Lemma~3.3]{InfiniteLanguages}.
Map the root to the root. If $f\in T_i(n+1)$, the three data
\[
\Psi_i(f\res n),\qquad \Psi_{i+1}(f_n),\qquad f(n)
\]
determine an immediate ambient successor of $\Psi_i(f\res n)$. Define
$\Psi_i(f)$ to be the unique next-level node of $S_i$ in that cone.
The recursion is well founded because both required map values have source
height $n$. For a strong vector subtree of height $h$ in the first $h$
coordinates, the same definition applies on $T_i(<h-i)$.

\begin{lemma}\label{lem:maps}
The maps $\Psi_i$ preserve meets and satisfy
\begin{align}
|\Psi_i(f)|&=e(|f|), &
\Psi_i(f)\res e(m)&=\Psi_i(f\res m),\label{eq:maps1}\\
(\Psi_i(f))_{e(m)}&=\Psi_{i+1}(f_m), &
\Psi_i(f)(e(m))&=f(m).\label{eq:maps2}
\end{align}
The restriction identity allows $m\leq|f|$; the slice and scalar identities
require $m<|f|$. For every finite $Y$ in the top domain,
\begin{equation}\label{eq:covariance}
C_i(\Psi_0[Y])=\Psi_i[C_i(Y)],\qquad
L(\Psi_0[Y])=e[L(Y)].
\end{equation}
Consequently $(u,f)\mapsto(u,\Psi_0(f))$ preserves all normal forms.
\end{lemma}

\begin{proof}
Induct on the source height, simultaneously for all $i$. The prescribed
successor cone gives~\eqref{eq:maps1}--\eqref{eq:maps2}; subsequent extensions
preserve the entire slice, including its values off the selected level set.
If $|f\wedge g|=c$ and $f,g$ are incomparable, then
\[
 f(c)\ne g(c)\quad\text{or}\quad f_c\ne g_c.
\]
In either case induction places $\Psi_i(f)$ and $\Psi_i(g)$ in different
immediate cones above $\Psi_i(f\res c)$. Therefore
\[
 \Psi_i(f)\wedge\Psi_i(g)=\Psi_i(f\wedge g).
\]
The generating operations commute with $\Psi$, so minimality and induction
on derivations give~\eqref{eq:covariance}. For finite vectors,~\eqref{eq:triangle}
places every argument in its stated domain. Finally, repeated use of
\eqref{eq:maps2} gives
\[
 \Psi_i(f)(e(a_1),\ldots,e(a_r))=f(a_1,\ldots,a_r),
\]
so compression, labels and normal forms are preserved.
\end{proof}

\subsection{Finite completion and the Ramsey property}

For a finite strong vector subtree $V=(V_i)_{i<h}$ of height $h$, the top
image $\Psi_0^V[T_0(<h)]$ is precisely $\val(V_0,\ldots,V_{h-1})$ in
\cite[Definition~3.3]{InfiniteLanguages}. Indeed, its successors are specified
by a node in the preceding top level, a node in the corresponding level of
the tail valuation tree, and a scalar. This is the defining recursion there.
We call this image a valuation tree.

The following relative completion uses the envelope reconstruction
of~\cite[Lemma~4.4]{InfiniteLanguages}, originating in the two-tree
argument of~\cite[Lemma~3.3]{TernaryFinite}. The closure and covariance
identities verify the hypotheses needed here.

\begin{lemma}[Relative finite completion]\label{lem:completion}
Let $\Psi$ arise from an infinite strong vector subtree $(S_i)_{i<\omega}$,
and let $X$ be a nonempty finite subset of its top image. Put $h=|L(X)|$.
There is $V\in\Str_h(S_0,\ldots,S_{h-1})$, also strong in the ambient
trees, with ambient level set $L(X)$, such that
\begin{equation}\label{eq:reconstruction}
\Psi_i^V(\widehat f)=f\qquad(f\in C_i(X)).
\end{equation}
In particular $X\subseteq\val(V)$. The conclusion also holds for an arbitrary
finite $X\subseteq T_0$, by taking $S_i=T_i$.
\end{lemma}

\begin{proof}
Write $L(X)=\{\ell_0<\cdots<\ell_{h-1}\}$. By~\eqref{eq:covariance},
$C_i(X)\subseteq S_i$. Every nonempty $C_i(X)$ has a unique node at
$\ell_0$: existence follows from restriction, and two such nodes would have
a meet below $\ell_0$.

For each $i<h$, complete $C_i(X)$ to $V_i$ on these $h$ levels. At an
unprescribed root choose any node. At consecutive selected levels
$\ell_j<\ell_{j+1}$ retain a prescribed successor when present and fill
every remaining immediate $S_i$-cone arbitrarily. There cannot be two
prescribed successors in one cone, since their meet would satisfy
\[
 \ell_j<|u\wedge v|<\ell_{j+1},\qquad |u\wedge v|\in L(X).
\]
Every prescribed descendant has its earlier restrictions in $C_i(X)$.
Thus the construction retains $C_i(X)$ and gives a strong subtree of $S_i$;
it is also strong in $T_i$. This is the completion argument of
\cite[Observation~2.1]{InfiniteLanguages}.

It remains to prove~\eqref{eq:reconstruction}. Induct on $j$, simultaneously
for $i+j<h$. The assertion is immediate at $j=0$. For
$f\in C_i(X)$ of height $\ell_j$, $j>0$, put
\[
 p=f\res\ell_{j-1},\qquad u=f_{\ell_{j-1}}.
\]
Then $p\in C_i(X)$, $u\in C_{i+1}(X)$, and induction gives
\[
 \Psi_i^V(\widehat p)=p,\qquad
 \Psi_{i+1}^V(\widehat u)=u,\qquad
 \widehat f(j-1)=f(\ell_{j-1}).
\]
These are precisely the predecessor, full slice and scalar specifying the
cone of $f$. Its unique next $V_i$-node is $f$, proving the induction step.
\end{proof}

\begin{proposition}\label{prop:source}
For every finite nonempty $U$, the structure $N_{\sigma,U}$ has the infinite Ramsey
property. Every finite member of its age is rigid.
\end{proposition}

\begin{proof}
We use the Milliken reduction of~\cite[proof of Proposition~4.1]{InfiniteLanguages}.

Let $B=\{b_0,\ldots,b_{k-1}\}\in\Age(N_{\sigma,U})$ be nonempty, and let
$c:\Emb(B,N_{\sigma,U})\to r$. Fix a realization of $B$ with normal form
$\tau$, compressed tuple $P_\tau=((u_j,f_j))_{j<k}$ and closure height $h$.
The relation naming $\tau$ determines every ordered subtuple type, so an
embedding of $B$ is precisely a realization of $\tau$. Define
\[
 \widetilde c(V)=c\bigl(b_j\mapsto(u_j,\Psi_0^V(f_j))\bigr),
 \qquad V\in\Str_h(T_0,\ldots,T_{h-1}).
\]
This is well defined by Lemma~\ref{lem:maps}.
Milliken's theorem gives $(S_i)_{i<h}$ on which $\widetilde c$ is constant.
Extend the vector to all coordinates with the same level set and let
\[
 z(u,f)=(u,\Psi_0(f)).
\]
Then $z\in\Emb(N_{\sigma,U})$. Let $a:B\to N_{\sigma,U}$ be an embedding.
By Lemma~\ref{lem:completion} there is
\[
 V\in\Str_h(S_0,\ldots,S_{h-1})
\]
such that
\[
 za(b_j)=(u_j,\Psi_0^V(f_j))\quad(j<k).
\]
Consequently $c(za)=\widetilde c(V)$ is constant.

Every finite member is rigid: the normal form of a tuple listing all its
vertices records an ordered list of pairwise distinct labelled compressed
nodes, by Lemma~\ref{lem:compression}. The empty case is immediate.
\end{proof}

\section{A uniform embedding into the valuation hypergraph}\label{sec:envelopes}

We now assume $\sigma_j=n_{j+1}+2\geq2$. In $G_\sigma$, a set
$\{f_0,\ldots,f_{r-1}\}$ of nodes of strictly decreasing heights has colour
$s\in\{1,\ldots,n_r\}$ exactly when
\[
f_0(|f_1|,\ldots,|f_{r-1}|)=s.
\]
The values $0$ and $\sigma_{r-1}-1$ both represent no relation, as in
\cite[Definition~3.1]{InfiniteLanguages}. Sets with height ties have no
relation. An envelope of a finite set is a valuation tree containing it,
with no requirement of minimality~\cite[Definition~4.3]{InfiniteLanguages}.
If the language has a largest arity, we keep the unused alphabet sizes
equal to $2$, whereas Section~4 of that paper sets them to $1$.
For unbounded arity, and in particular under the hypothesis of
Conjecture~4.7, our signature agrees exactly with theirs.

Enumerate $H$ by $\N$. For $v>z_0>\cdots>z_{t-1}$ let
$c(v,\bar z)\in\{0,\ldots,\sigma_t-2\}$ be its genuine colour.
For every nonempty decreasing tuple $\bar y$ put
\begin{equation}\label{eq:cutoff}
D(\bar y)=\max_{0\leq d\leq y_0}(\sigma_{d+|\bar y|}-2).
\end{equation}
Introduce distinct markers $b_s(\bar y)$ for $1\leq s\leq D(\bar y)$.
Place those beginning with $n$ in a finite block immediately before vertex
$n$, with an arbitrary fixed order inside each block. Every block is finite,
so the resulting order has type $\omega$.
Identify its positions with $\N$, denote original positions by $o_v$, and put
$O=\{o_v:v<\omega\}$ and $B=\N\setminus O$. In particular
$b_s(\bar y)<o_{y_0}$.
An $m$-slice means successive slicing along a decreasing prefix of length
$m$; the $0$-slice is the node itself. A slice is \emph{original} if every
position in its prefix belongs to $O$.

Define $\phi(v)\in T_0(o_v)$ by
\begin{align}
\phi(v)(o_{z_0},\ldots,o_{z_{t-1}})&=c(v,\bar z),\label{eq:original-code}\\
\phi(v)(o_{x_0},\ldots,o_{x_{m-1}},b_s(\bar y))
 &=\sigma_{m+1}-1\label{eq:marker-code}
\end{align}
in the second line precisely when $(v)\cat\bar x\cat\bar y$ is strictly
decreasing, $m=|\bar x|\leq y_0$, and $c(v,\bar x\cat\bar y)=s>0$.
The prefix $\bar x$ may be empty. All other entries are zero.
Equation~\eqref{eq:cutoff} ensures that every required marker exists.
The prescribed value in~\eqref{eq:marker-code} lies in the correct alphabet.
Equation~\eqref{eq:original-code} proves that $\phi$ is an embedding, and its
strictly increasing heights preserve the enumeration.

This modifies the fixed-depth marker construction
of~\cite[Lemma~4.3]{InfiniteLanguages}. The cutoff now depends on the first
vertex of the suffix, rather than on a chosen size of finite substructure.

\begin{lemma}\label{lem:weak-meets}
Every slice of an image node whose prefix contains a position in $B$ is zero.
Suppose that $f,g\in T_m$ are incomparable and each is either an original
$m$-slice of an image node or a zero node. If $|f\wedge g|=o_n$, then $n<m$.
\end{lemma}

\begin{proof}
A nonzero entry containing a marker has exactly one marker, in its last
position. If a slice prefix contains a marker, every entry of that slice
therefore vanishes.

Let $f,g\in T_m$ be the nodes in the second assertion; allow one to be
zero. Suppose, towards a contradiction, that
\[
 |f\wedge g|=o_n,\qquad n\geq m.
\]
Then $f\res o_n=g\res o_n$, but there is a distinguishing input whose
largest coordinate is $o_n$. Each nonzero comparator has an original outer
prefix of length $m$, all of whose vertices exceed $n$.

\emph{Case 1.} The input is $(o_{z_0},\ldots,o_{z_{t-1}})$, where $z_0=n$.
Choose a genuine colour $s>0$ taken on precisely one side. Since
$m\leq n$,~\eqref{eq:cutoff} supplies the marker $b=b_s(\bar z)$.

\emph{Case 2.} The input contains a marker. A nonzero value requires the form
\[
 (o_{z_0},\ldots,o_{z_{q-1}},b_s(\bar y)),\qquad
 q\geq1,\quad z_0=n,\quad m+q\leq y_0.
\]
The tuple $\bar z\cat\bar y$ is decreasing. For two original comparators
these admissibility conditions coincide; hence exactly one corresponding
genuine colour with suffix $\bar z\cat\bar y$ equals $s$. The same conclusion
holds against zero. Since $m\leq n$, the marker
$b=b_s(\bar z\cat\bar y)$ exists.

In either case $b<o_n$, and~\eqref{eq:marker-code} gives
\[
 \{f(b),g(b)\}=\{0,\sigma_{m+1}-1\},
\]
contrary to $f\res o_n=g\res o_n$.
These cases exhaust the possible nonzero entries.
\end{proof}

\begin{lemma}\label{lem:bound}
For every $S\in[H]^k$, $k\geq1$, the intrinsic closure of $\phi[S]$ has
at most $R(k)$ levels.
\end{lemma}

\begin{proof}
The finite container below replaces the closure construction
of~\cite[Lemma~4.4]{InfiniteLanguages}.
Put $A=S\cup\{0,\ldots,k-1\}$ and $N=\max o[S]$. Then $|A|\leq2k$
and $\max A=\max S$. Let $F_i$ be the family of original slices
\[
\phi(v)_{o_{a_1},\ldots,o_{a_i}},\qquad
v\in S,\quad v>a_1>\cdots>a_i,\quad a_j\in A.
\]
Thus $F_0=\phi[S]$, $F_i=\varnothing$ for $i\geq|A|$, and
\begin{equation}\label{eq:P}
P=\sum_i|F_i|\leq k2^{|A|}\leq k4^k.
\end{equation}
Write $0_N^i$ for the zero node at height $N$ in coordinate $i$ and set
\begin{equation}\label{eq:Lstar}
L^*=o[A]\cup\bigcup_{i:F_i\ne\varnothing}
\{|u\wedge w|:u,w\in F_i\cup\{0_N^i\}\}.
\end{equation}

We claim that
\begin{equation}\label{eq:original-level-trap}
 L^*\cap O\subseteq o[A].
\end{equation}
A comparable meet has one of the original heights already in $o[A]$.
For an incomparable meet $|u\wedge w|=o_n$, there are two cases.
If either original prefix uses $a<k$, then
$o_n<\min\{|u|,|w|\}\leq o_a<o_k$, so $n<k$.
Otherwise the prefixes use only $S$, whence $i<k$ and
Lemma~\ref{lem:weak-meets} gives $n<i<k$. These statements include the
comparison with $0_N^i$, proving~\eqref{eq:original-level-trap}.

Put $p_i=|F_i|$. The zero--zero meet contributes $N\in o[A]$, and thus
\begin{equation}\label{eq:level-count}
|L^*|\leq2k+\sum_i(p_i^2+2p_i)\leq2k+P^2+2P\leq R(k).
\end{equation}

Let $L^*=\{\ell_0<\cdots<\ell_{h-1}\}$ and define, for $i<h$,
\begin{equation}\label{eq:container}
Q_i=\{u\res\ell_j:u\in F_i,\ \ell_j\leq|u|\}
     \cup\{0_{\ell_j}^i:i+j<h\}.
\end{equation}
Every $u\in F_i$ has $i$ higher levels in $o[A]$, namely its original
vertex and all but the last slicing parameter. Consequently
\begin{equation}\label{eq:container-triangle}
 f\in Q_i,\quad |f|=\ell_j\quad\Longrightarrow\quad i+j<h.
\end{equation}
We verify closure of $Q$.
\begin{enumerate}
\item Restrictions remain in~\eqref{eq:container}.
\item For $u,w\in F_i\cup\{0_N^i\}$ and admissible $a,b\in L^*$,
\[
 (u\res a)\wedge(w\res b)
   =u\res\min\{a,b,|u\wedge w|\}\in Q_i.
\]
The minimum belongs to $L^*$ by~\eqref{eq:Lstar}.
\item Let $f\in Q_i$ have height $\ell_j$, and let $t<j$.
If $f$ is zero or $\ell_t\in B$, then
\[
 f_{\ell_t}=0^{{i+1}}_{\ell_t}\in Q_{i+1},
 \qquad i+1+t\leq i+j<h.
\]
Otherwise $f$ is a restriction of some $u\in F_i$ and
$\ell_t=o_a$ for $a\in A$, by~\eqref{eq:original-level-trap}. Then
$f_{\ell_t}=u_{o_a}\in F_{i+1}\subseteq Q_{i+1}$.
\end{enumerate}
Thus $\phi[S]\subseteq Q_0$ and $C_i(\phi[S])\subseteq Q_i$ for every $i$.
In particular $L(\phi[S])\subseteq L^*$, and~\eqref{eq:level-count} finishes
the proof.
\end{proof}

\begin{proof}[Proof of Theorem~\ref{thm:envelopes} and Corollary~\ref{cor:counterexample}]
Apply Lemma~\ref{lem:completion} in the full ambient trees to $\phi[S]$.
It gives a genuine valuation envelope of height $|L(\phi[S])|\leq R(k)$.
The map $\phi$ was chosen independently of $k$ and $S$.
If $n_r$ is unbounded and $H$ is universal, choosing a threshold greater than
$R(k)$ contradicts Conjecture~4.7 for every positive $k$; for $k=0$, any
valuation tree of height one contains the empty set.
\end{proof}

\begin{remark}
The stronger $k$-enveloping condition of
\cite[Definition~4.2]{InfiniteLanguages} excludes every original meet of
incomparable slices. Lemma~\ref{lem:weak-meets} permits $n<m$; by
\eqref{eq:original-level-trap}, these exceptions do not prevent uniformly
bounded envelopes.
\end{remark}

\section{From hypergraphs to relational structures}\label{sec:profiles}

Return to the language of Corollary~\ref{thm:brs}. Write $L_r$ for its
$r$-ary symbols and let $U=\{0,1\}^{L_1}$. For $r\geq2$, on an increasingly
enumerated $r$-element support record the full profile
\[
\bigl(R(a_{\pi(0)},\ldots,a_{\pi(r-1)})\bigr)_{R\in L_r,\,\pi\in\Sym_r}.
\]
There are $q_r=2^{|L_r|r!}$ profiles. Identify them with $\{0,\ldots,q_r-1\}$,
taking the all-false profile to zero, and put $\sigma_{r-1}=q_r+1$.
This is the usual reduction to coloured hypergraphs
used in~\cite[Proposition~5.2]{InfiniteLanguages}.

Define an $L$-structure $M$ on $U\times T_0$. Unary relations read the
label. For an $r$-element support of distinct heights, list the vertices in
increasing height order and evaluate the highest node on the other heights
in decreasing order. Decode a value below $q_r$ as its profile and decode
the reserved value $q_r$ as the all-false profile. On supports with a height
tie, make every relation of arity at least two false. Relations on repeated
entries are always false.

An enumeration of $K$ makes its support profiles a coloured hypergraph.
Apply the construction of Section~\ref{sec:envelopes} and retain the unary
labels. This gives an embedding $\phi:K\to M$. It preserves the chosen
enumeration by heights, so every orientation in a profile is decoded in
the correct order. Every finite substructure of $M$ is an injective $L$-structure, so
$\Age(M)\subseteq\Age(K)$. Countable universality of the Fra\"iss\'e
limit gives an embedding $\theta:M\to K$.

Expand $M$ by the normal-form predicates of $N_{\sigma,U}$, and call the
result $N$. Every $L$-relation is determined by the normal form of its
tuple: that form records labels, height order and ties, and the valuation
on the support. Thus adding the $L$-relations does not change embeddings
between the relevant finite expanded structures or self-embeddings of the
source. Proposition~\ref{prop:source} applies to $N$.

\begin{proposition}\label{prop:precompact}
The expansion $X=N\cdot\phi$ of $K$ is precompact. For every $k$, only
finitely many of its added $k$-ary predicates are nonempty. If $A$ has
$k\geq1$ vertices and
\[
B_\sigma(R)=\sum_{n<R}\prod_{r=1}^n\sigma_r^{\binom nr},
\]
then
\begin{equation}\label{eq:crude-bound}
|\mathcal P_A(X)|\leq\bigl(|U|B_\sigma(R(k))\bigr)^k.
\end{equation}
\end{proposition}

\begin{proof}
Lemma~\ref{lem:bound} bounds the closure height of every projected
$k$-tuple in $\phi[K]$. Compression therefore places its projected normal
form in the finite set $T_0(<R(k))^k$. Including labels gives the displayed
bound. By heredity of normal forms, the whole pullback expansion on the
tuple is determined by this one normal form. Only finitely many such
forms occur globally at each size, which also proves the assertion about
nonempty predicates.
The empty structure has a unique pullback type.
\end{proof}
\section{Compatible optimal expansions}\label{sec:extraction}

We prove Theorem~\ref{thm:transfer} using the upper-bound transfer
\cite[Theorem~3.10]{HZ}, finite-palette minimization
\cite[Fact~3.4]{HZ}, and the recurrence argument in the proof of
\cite[Theorem~3.18]{HZ}. Complete source types are first named by finitely
many predicates at each stage. Recurrence then makes the stages extend
one another on the same underlying set. We use only the recurrence
argument of the last reference; optimality below is proved separately on
the finite family handled at each stage.

Fix $K,M,N,\phi,\theta$ as in Theorem~\ref{thm:transfer}, and put
$X=N\cdot\phi$. Choose one representative of each isomorphism type of
nonempty finite substructures of $K$, labelling each $k$-point representative
by $\{0,\ldots,k-1\}$. Take increasing
finite sets $I_0=\varnothing\subseteq I_1\subseteq\cdots$ exhausting
them and closed under nonempty substructure types. Such an exhaustion
exists because $K$ is countable and a finite structure has only finitely
many induced substructures.

For each representative $A$, put $\mathcal Q_A=\mathcal P_A(X)$.
For $p\in\mathcal Q_A$, introduce an injective $|A|$-ary predicate
$R_{A,p}$. Let $N_n$ expand $M$ by the symbols indexed by $A\in I_n$.
Writing $k=|A|$, set
\[
R_{A,p}^{N_n}(f(0),\ldots,f(k-1))
\quad\Longleftrightarrow\quad
f\in\Emb(A,M)\ \text{ and }\ N\cdot f=p.
\]
All other tuples give false instances. Set $Y_n=N_n\cdot\phi$, so
$Y_0=K$. Every $N$-self-embedding preserves these predicates; hence
$\Emb(N)\subseteq\Emb(N_n)$.

For $\eta\in\Emb(K)$, put $F=Y_n\cdot\eta$. If $A\in I_n$ and
$a,b\in\Emb(A,K)$, then
\begin{equation}\label{eq:source-type-identification}
F\cdot a=F\cdot b
\quad\Longleftrightarrow\quad
X\cdot(\eta a)=X\cdot(\eta b).
\end{equation}
Indeed, $R_{A,p}$ on the fixed enumeration records the complete source
type, which in turn determines all predicates on ordered subtuples.
Consequently each $p\in\mathcal P_A(F)$ determines a unique
$A^p\in\mathcal Q_A$, and
\begin{equation}\label{eq:recoded-source}
\Emb(p,N_n)=\Emb(A^p,N),\qquad \brd(p,N_n)=1.
\end{equation}
The equality follows from the distinguished positive naming predicate.
The degree-one assertion follows from the hypothesis on $A^p$ and
$\Emb(N)\subseteq\Emb(N_n)$.

\begin{lemma}\label{lem:finite-stage}
Let $n\geq0$, $\eta\in\Emb(K)$ and $F=Y_n\cdot\eta$.
If $\mathcal P_A(F)$ is unavoidable for every $A\in I_n$, then $F$ is
recurrent and $|\mathcal P_A(F)|=\brd(A,K)$ for every $A\in I_n$.
For a fixed $A\in I_n$, the upper bound does not require unavoidability.
\end{lemma}

\begin{proof}
By~\eqref{eq:recoded-source}, the upper-bound transfer
\cite[Theorem~3.10]{HZ}, applied to $(M,\phi\eta,\theta)$ and $N_n$,
gives
\[
\brd(A,K)\leq
\sum_{p\in\mathcal P_A(F)}\brd(p,N_n)=|\mathcal P_A(F)|.
\]
Unavoidability gives the reverse inequality by~\cite[Fact~3.4(1)]{HZ}.

For recurrence we use the argument of~\cite[proof of Theorem~3.18]{HZ}.
Fix $g\in\Emb(K)$. Simultaneously homogenize the finitely many colourings
\[
d_{A,p}:\Emb(p,N_n)\longrightarrow\mathcal P_A(F),\qquad
d_{A,p}(f)=F\cdot(g\theta f),
\quad A\in I_n, p\in\mathcal P_A(F),
\]
using~\eqref{eq:recoded-source}. Let $z\in\Emb(N_n)$ witness this and
put $e=g\theta z\phi\eta$. There are maps $\pi_A$ satisfying
\begin{equation}\label{eq:palette-permutation}
F\cdot(ea)=\pi_A(F\cdot a),\qquad
\pi_A[\mathcal P_A(F)]=\mathcal P_A(F\cdot e)=\mathcal P_A(F).
\end{equation}
Thus the $\pi_A$ are permutations. Choose a common multiple $q\geq1$
of their orders. Then $F\cdot(e^qa)=F\cdot a$ for all $A\in I_n$.
These types determine every added predicate; on repeated tuples and
tuples of the wrong base type the predicates are false. Therefore
$e^q\in\Emb(F)$. Since
$e^q=g[(\theta z\phi\eta)e^{q-1}]$, recurrence follows. The case
$I_n=\varnothing$ is immediate.
\end{proof}

\begin{proof}[Proof of Theorem~\ref{thm:transfer}]
Construct $\eta_n\in\Emb(K)$ and $F_n=Y_n\cdot\eta_n$ so that
\begin{enumerate}
\item $F_n$ is recurrent;
\item its palettes on $I_n$ are optimal and unavoidable;
\item $F_n$ extends $F_{n-1}$ on the same underlying set.
\end{enumerate}
Start with $\eta_0=\mathrm{id}_K$ and $F_0=K$. Given the preceding
stage, put $H_n=Y_n\cdot\eta_{n-1}$ and choose $g\in\Emb(K)$
minimizing
\begin{equation}\label{eq:minsum}
\sum_{A\in I_n}|\mathcal P_A(H_n\cdot g)|.
\end{equation}
For every $u\in\Emb(K)$, the inclusions
$\mathcal P_A(H_n\cdot gu)\subseteq\mathcal P_A(H_n\cdot g)$
are equalities by minimality. Hence these palettes are unavoidable;
this is~\cite[Fact~3.4(2)]{HZ}.

By recurrence of $F_{n-1}$ choose $h\in\Emb(K)$ with
$gh\in\Emb(F_{n-1})$, and set
\[
\eta_n=\eta_{n-1}gh,\qquad F_n=Y_n\cdot\eta_n.
\]
The palettes remain unchanged, so Lemma~\ref{lem:finite-stage} gives
(1)--(2). Moreover,
\begin{equation}\label{eq:literal}
F_n\res L_{n-1}^*=F_{n-1}\cdot gh=F_{n-1},
\end{equation}
where $L_{n-1}^*$ is the language of $F_{n-1}$, proving (3).

Let $K^*=\bigcup_nF_n$. If $A\in I_n$, a later predicate $R_{B,p}$
has $B\notin I_n$. Such a $B$ cannot embed into $A$ by closure of
$I_n$ under substructure types. Every later predicate is therefore
false on an $A$-copy, and~\eqref{eq:literal} yields
\begin{equation}\label{eq:final-type-identity}
K^*\cdot a=F_n\cdot a\qquad(a\in\Emb(A,K)),
\end{equation}
after restoring the empty later predicates. Thus the final palettes
are finite, optimal and unavoidable, so $K^*$ is a big Ramsey structure.
The empty structure has degree one.

By~\eqref{eq:source-type-identification}, the palette of
$Y_n\cdot\eta$ on $A\in I_n$ is in bijection with
$\mathcal P_A(X\cdot\eta)$. Lemma~\ref{lem:finite-stage} gives the
upper bound in~\eqref{eq:transfer-min} for every $\eta$, and equality
holds at $\eta=\eta_n$ by optimality. This proves the minimum formula.
Finally, the pairs $(A,p)$ form a countable set. If there are finitely
many base types of each finite size, only finitely many such pairs
have $|A|=k$, giving finitely many added $k$-ary symbols.
\end{proof}

The union above is a union of compatible interpretations on $K$;
no intersection of the ranges of $\eta_n$ is taken. We do not assert
recurrence of $K^*$ or simultaneous optimality of one pullback
$X\cdot\eta$ for all finite structures.

\begin{proof}[Proof of Corollary~\ref{thm:brs}]
Use $M,N,\phi,\theta$ from Section~\ref{sec:profiles}.
Proposition~\ref{prop:precompact} gives precompactness of
$X=N\cdot\phi$, while Proposition~\ref{prop:source} gives degree one
for every finite member of $\Age(N)$, hence of $\Age(X)$.
Theorem~\ref{thm:transfer} applies. For each $k$, only the finitely
many symbols of arity at most $k$ can be nonempty on an injective
$k$-point structure. Thus there are finitely many base types of each
size, and the theorem's language-size conclusion applies.
\end{proof}

Return now to the structures and embeddings of Section~\ref{sec:profiles}.

\begin{proposition}\label{prop:degrees}
For finite nonempty $A\in\Age(K)$ with vertices labelled
$0,\ldots,k-1$, put
$\chi_A(e)=\NF(\phi e(0),\ldots,\phi e(k-1))$. Then
\begin{equation}\label{eq:exact}
\brd(A,K)=\min_{\eta\in\Emb(K)}
|\chi_A[\eta\circ\Emb(A,K)]|.
\end{equation}
Every finite colouring of $\Emb(A,K)$ has a self-copy on which its
value depends only on the $K^*$-type. Moreover,
\begin{equation}\label{eq:copy}
\crd(A,K)=\frac{\brd(A,K)}{|\Aut(A)|}.
\end{equation}
\end{proposition}

\begin{proof}
Normal forms determine precisely the full pullback types in $X$, so
\eqref{eq:exact} is~\eqref{eq:transfer-min}. For canonization, put
$d=|\mathcal P_A(K^*)|$ and apply the degree bound to
$e\mapsto(c(e),K^*\cdot e)$. On a self-copy there are at most $d$
pairs, and unavoidability forces all $d$ types to occur. Each type
therefore has exactly one $c$-colour; compare~\cite[Proposition~3.8]{HZ}.
Finally,~\eqref{eq:copy} is the standard
identity~\eqref{eq:degree-conventions}.
\end{proof}
\section{A finitely constrained class with infinite singleton degree}
\label{neg:section}

The unrestricted classes considered above have finite big Ramsey degrees.
Finite forbidden families produce a different situation. We construct one
finite language in which constraints on at most four vertices force infinite
degree already for a singleton. The obstruction is expressed by the ages of
one-point realization sets, rather than by the sizes of coding envelopes.
These ages are the ranks of types in Sauer's work on vertex partitions
\cite{SauerHomogeneous2020}. The passage from incomparable ages to
persistent vertex colours is Sauer's persistent-partition argument
\cite[Section~5]{SauerCanonicalVertex2003}. Here the construction realizes
an infinite antichain of these ages in one finitely constrained class.

A finite relational structure is \emph{irreducible} if each pair of distinct
vertices occurs together in a positive relation tuple. Write
$\operatorname{Forb}_{e}(\mathcal F)$ for the finite structures admitting no
embedding of a member of $\mathcal F$. In this section, embeddings
preserve and reflect relations, and relations hold only on tuples of
distinct vertices.

The following form also gives the conclusion for larger finite structures
containing the distinguished singleton.

\begin{theorem}\label{neg:main}
There are a finite relational language $L_*$ of maximum arity three, a finite
family $\mathcal F_*$ of irreducible $L_*$-structures, each with three or four
vertices, and a one-point structure $\delta$ such that the following holds.
The class $\operatorname{Forb}_{e}(\mathcal F_*)$ is a free amalgamation
class with a countably infinite Fra\"iss\'e limit $K_*$, and
\[
 \brd(\delta,K_*)=\crd(\delta,K_*)=\infty.
\]
More generally, both degrees are infinite for every finite
$A\in\Age(K_*)$ containing a vertex of type $\delta$.
\end{theorem}

The language and the forbidden family in this theorem are fixed. In
particular, they do not depend on the number of colours used to witness
infinite degree. The theorem supplies an obstruction relevant to
Problem~6.3.3 of~\cite{BinaryExact}; it does not classify the finitely
constrained free amalgamation classes.

\subsection{Realization sets and persistent colourings}
\label{neg:types}

Let $M$ be a countable homogeneous structure in a finite relational language
whose age has free amalgamation. A \emph{one-point type} over a finite
$A\subseteq M$ is a structure $t$ on $A\cup\{*\}$ extending $M\res A$
and belonging to $\Age(M)$. Its realization set and its age are
\begin{align*}
 P(t)&=\{x\in M\setminus A:
       \operatorname{id}_{A}\cup\{(*,x)\}\in\Emb(t,M)\},\\
 \mathcal C(t)&=\Age(M\res P(t)).
\end{align*}
We call a structure \emph{age indivisible} if, in every finite partition of
its vertex set, some part has the same age as the whole structure.

The following age-indivisibility result applies to types in arbitrary
finite relational languages.

\begin{lemma}\label{neg:age-indivisibility}
For every realized one-point type $t$, the structure $M\res P(t)$ is age
indivisible. The same holds for every countable structure with age
$\mathcal C(t)$.
\end{lemma}

\begin{proof}
Homogeneity makes $P(t)$ an orbit of the pointwise stabilizer of its
finite parameter set. By~\cite[Theorem~3.2 and Lemma~4.1]{SauerCanonicalVertex2003},
free amalgamation gives age indivisibility of this orbit.
The property depends only on the age by
\cite[Theorem~3.1]{SauerCanonicalVertex2003}, proving the second assertion.
The finite-colour statement follows from the two-colour statement by induction.
\end{proof}

The next proposition is a formulation of Sauer's persistence argument
\cite[Lemma~5.1 and Theorem~5.1]{SauerCanonicalVertex2003}. We retain the
short proof for the specified singleton type $\delta$.

\begin{proposition}\label{neg:persistence}
For $1\leq i\leq q$, let $t_i$ be a realized one-point type over a finite
set $A_i\subseteq M$. Suppose their points have the same singleton
type $\delta$ and their ages $\mathcal C_i=\mathcal C(t_i)$ are pairwise
incomparable under inclusion. There is a colouring
\[
 \chi:M_\delta\longrightarrow\{1,\ldots,q\},\qquad
 M_\delta=\{x\in M:M\res\{x\}\cong\delta\},
\]
such that, for every $f\in\Emb(M,M)$ and $i\le q$,
\begin{equation}\label{neg:persistent-age}
 \Age\bigl(M\res(f[P(t_i)]\cap\chi^{-1}(i))\bigr)=\mathcal C_i.
\end{equation}
In particular, all $q$ colours occur in every copy of $M$.
\end{proposition}

\begin{proof}
Fix an enumeration $\nu:M\to\omega$ and, for finite $B\subseteq M$, put
\[
 w(B)=\sum_{x\in B}2^{\nu(x)}.
\]
For $x\in M_\delta$, minimize $w(\alpha[A_i])$ over all
$i\le q$ and all $\alpha\in\Emb(t_i,M)$ with $\alpha(*)=x$.
Such embeddings exist by homogeneity. If two indices attain the same
minimum, their boundary images coincide, since $w$ is injective. The two
pointed types are then isomorphic. Homogeneity implies that their
realization sets have equal ages, contrary to incomparability. Let
$\chi(x)$ be the unique minimizing index.

Fix $i$ and $f$, and write
$Y=f[P(t_i)]$, $B=f[A_i]$. If $x\in Y$ and $\chi(x)\ne i$, a minimizing
witness has index $j\ne i$ and boundary of weight strictly less than
$w(B)$. Only finitely many transported realization sets
$Q_1,\ldots,Q_m$ can arise from these witnesses. For each $a\le m$,
\[
 \Age(M\res Q_a)=\mathcal C_{j_a},\qquad j_a\ne i,
 \qquad D_a\in\mathcal C_i\setminus\mathcal C_{j_a}
\]
for some finite $D_a$. Partition $Y$ into
\[
 Z_a=(Y\cap Q_a)\setminus\bigcup_{b<a}Q_b\quad(1\le a\le m),
 \qquad Z_0=Y\setminus\bigcup_{a=1}^{m}Q_a.
\]
No $Z_a$ with $a\ge1$ contains $D_a$. Since $\Age(M\res Y)=\mathcal C_i$,
Lemma~\ref{neg:age-indivisibility} yields
$\Age(M\res Z_0)=\mathcal C_i$. By construction $\chi[Z_0]=\{i\}$,
which proves~\eqref{neg:persistent-age}.
\end{proof}

\subsection{The forbidden family}
\label{neg:family}

We now encode the homomorphism order of finite directed graphs in the ages
of realization sets. The language is
\[
\begin{array}{c|l}
\text{arity}&\text{symbols}\\ \hline
1&O,V,\mathsf{Path},\mathsf{Test},\mathsf{First},\mathsf{Last}\\
2&E,N,Q,G,L,R\\
3&S_0,S_1,S_2,U_L,U_R.
\end{array}
\]
An \emph{$E$-complete} set satisfies $E(x,y)$ for all distinct $x,y$ in it.
Put
\[
 \mathcal T=\{(0,0),(0,1),(1,2),(2,2)\}.
\]
The following conditions describe the legal structures. Every displayed
tuple of named vertices consists of distinct vertices and is $E$-complete.
An implication is required only when its displayed hypotheses hold.
\begin{enumerate}
\item On $\{a,h,v\}$, if
\[
 O(a)\wedge\mathsf{Path}(a)\wedge V(h)\wedge V(v)\wedge Q(h,v),
\]
then
\begin{equation}\label{neg:states}
 \sum_{s=0}^{2}\mathbf1_{S_s(a,h,v)}=1,\qquad
 \mathsf{First}(a)\Rightarrow S_0(a,h,v),\qquad
 \mathsf{Last}(a)\Rightarrow S_2(a,h,v).
\end{equation}
\item On $\{a,a',h,v\}$, if
\[
 O(a)\wedge\mathsf{Path}(a)\wedge O(a')\wedge\mathsf{Path}(a')
 \wedge V(h)\wedge V(v)\wedge N(a,a')\wedge Q(h,v),
\]
then
\begin{equation}\label{neg:transition}
 S_s(a,h,v)\wedge S_t(a',h,v)\ \Rightarrow\ (s,t)\in\mathcal T.
\end{equation}
\item On $\{b,a,h,v\}$, if
\[
 O(b)\wedge\mathsf{Test}(b)\wedge O(a)\wedge\mathsf{Path}(a)
 \wedge V(h)\wedge V(v)\wedge Q(h,v),
\]
then
\begin{equation}\label{neg:copy-bits}
 \begin{aligned}
 L(b,a)&\Rightarrow\bigl(U_L(b,h,v)\leftrightarrow S_1(a,h,v)\bigr),\\
 R(b,a)&\Rightarrow\bigl(U_R(b,h,v)\leftrightarrow S_1(a,h,v)\bigr).
 \end{aligned}
\end{equation}
\item On $\{b,h,v,w\}$, if
\[
 O(b)\wedge\mathsf{Test}(b)\wedge V(h)\wedge V(v)\wedge V(w)
 \wedge Q(h,v)\wedge Q(h,w)\wedge G(v,w),
\]
then
\begin{equation}\label{neg:conflict}
 \neg\bigl(U_L(b,h,v)\wedge U_R(b,h,w)\bigr).
\end{equation}
\end{enumerate}

For each condition, include in $\mathcal F_*$ every $E$-complete injective
$L_*$-structure on its displayed three- or four-element labelled domain
that violates that condition. Subject to $E$-completeness and the violation,
all other relation values are arbitrary.
Thus legality is exactly the satisfaction of
\eqref{neg:states}--\eqref{neg:conflict}. No additional assumptions about
the predicates or the binary and ternary relations are imposed.

\begin{lemma}\label{neg:fraisse}
The family $\mathcal F_*$ is finite and consists of irreducible structures
of size three or four. The class
$\mathcal K_*:=\operatorname{Forb}_{e}(\mathcal F_*)$ has a countably
infinite homogeneous Fra\"iss\'e limit $K_*$.
\end{lemma}

\begin{proof}
The language and the listed domains are finite. Every forbidden structure
is $E$-complete, hence irreducible. Free amalgamation follows
from~\cite[Theorem~3.2]{SauerCanonicalVertex2003}. The class is hereditary, has countably many
isomorphism types, and has arbitrarily large members with every relation
empty. Fra\"iss\'e's theorem~\cite{Fraisse} gives $K_*$.
\end{proof}

\subsection{Encoding directed-graph homomorphisms}
\label{neg:encoding}

Let $H$ be a nonempty finite loopless directed graph on $\{1,\ldots,k\}$.
Opposite arcs are permitted. A homomorphism preserves arcs and need not be
injective. Write $J\to H$ when such a map exists, and put
\[
 Z_H=\{(i,j)\in\{1,\ldots,k\}^2:(i,j)\notin E(H)\}.
\]
In particular, $(i,i)\in Z_H$ for every $i$.

Define the boundary $A_H$ on the disjoint union
\[
 \{a_0,\ldots,a_{k+1}\}\ \dot\cup\
 \{b_{ij}:(i,j)\in Z_H\}
\]
by the following tables:
\begin{align*}
 O^{A_H}&=A_H,&V^{A_H}&=\varnothing,\\
 \mathsf{Path}^{A_H}&=\{a_0,\ldots,a_{k+1}\},&
 \mathsf{Test}^{A_H}&=\{b_{ij}:(i,j)\in Z_H\},\\
 \mathsf{First}^{A_H}&=\{a_0\},&
 \mathsf{Last}^{A_H}&=\{a_{k+1}\},\\
 E^{A_H}&=\{(x,y):x,y\in A_H,\ x\ne y\},&
 N^{A_H}&=\{(a_l,a_{l+1}):0\le l\le k\},\\
 L^{A_H}&=\{(b_{ij},a_i):(i,j)\in Z_H\},&
 R^{A_H}&=\{(b_{ij},a_j):(i,j)\in Z_H\}.
\end{align*}
All other relation values are false. Let $\delta$ be the singleton on which
$V$ alone holds. The type $t_H$ over $A_H$ gives $*$ singleton type $\delta$,
sets $E(*,a)=E(a,*)=1$ for every $a\in A_H$, and makes every other mixed
relation value false. Both $A_H$ and $t_H$ are legal, since each forbidden
condition requires at least two $V$-vertices. Embed $A_H$ into $K_*$ and put
\[
 P_H=P(t_H),\qquad \mathcal C_H=\Age(K_*\res P_H).
\]

For a finite $L_*$-structure $B$ and $h\in B$, define a loopless directed
graph $J_h(B)$ by
\begin{align}
 D_h(B)&=\{v\in B\setminus\{h\}:Q(h,v)\wedge E(h,v)\wedge E(v,h)\},
 \label{neg:local-domain}\\
 V(J_h(B))&=D_h(B),\notag\\
 E(J_h(B))&=\{(v,w)\in D_h(B)^2:v\ne w,
            \ E(v,w)\wedge E(w,v)\wedge G(v,w)\}.
 \label{neg:local-arcs}
\end{align}
The empty graph has a homomorphism into every $H$.

\begin{proposition}\label{neg:age-characterization}
For every finite $L_*$-structure $B$,
\begin{equation}\label{neg:age-formula}
 B\in\mathcal C_H
 \quad\Longleftrightarrow\quad
 \bigl(\forall x\in B\; B\res\{x\}\cong\delta\bigr)
 \ \text{and}\ \bigl(\forall h\in B\;J_h(B)\to H\bigr).
\end{equation}
No $E$-completeness assumption on $B$ is required.
\end{proposition}

\begin{proof}
The extension property of $K_*$ gives the following criterion:
\begin{equation}\label{neg:completion}
 \begin{split}
 B\in\mathcal C_H\quad\Longleftrightarrow\quad
 &\text{there is }Y\in\mathcal K_*\text{ on }A_H\dot\cup B\text{ with}\\
 &Y\res B=B,\qquad
 Y\res(A_H\cup\{x\})=t_H\quad(x\in B),
 \end{split}
\end{equation}
where the last equality sends $*$ to $x$.

Assume such a $Y$ exists. Fix $h\in B$ and $v\in D_h(B)$. By
\eqref{neg:states} there is a unique $c_l(h,v)\in\{0,1,2\}$ such that
$S_{c_l(h,v)}(a_l,h,v)$ holds. Equations
\eqref{neg:states} and~\eqref{neg:transition} give
\[
 c_0(h,v)=0,\qquad c_{k+1}(h,v)=2,\qquad
 (c_l(h,v),c_{l+1}(h,v))\in\mathcal T.
\]
Consequently there is a unique $f_h(v)\in\{1,\ldots,k\}$ with
\begin{equation}\label{neg:marker}
 c_l(h,v)=
 \begin{cases}
 0&l<f_h(v),\\1&l=f_h(v),\\2&l>f_h(v).
 \end{cases}
\end{equation}
For every $(i,j)\in Z_H$, equation~\eqref{neg:copy-bits} yields
\begin{equation}\label{neg:marker-bits}
 U_L(b_{ij},h,v)\ \Longleftrightarrow\ f_h(v)=i,
 \qquad
 U_R(b_{ij},h,v)\ \Longleftrightarrow\ f_h(v)=j.
\end{equation}
All these applications are on $E$-complete sets: the boundary is
$E$-complete, every old--new pair has both $E$-directions, and
$v\in D_h(B)$ supplies the remaining pair.

If $(v,w)\in E(J_h(B))$ and
$(f_h(v),f_h(w))\in Z_H$, take $b=b_{f_h(v),f_h(w)}$.
The four points $b,h,v,w$ satisfy the hypotheses of
\eqref{neg:conflict}, whereas both prohibited bits hold by
\eqref{neg:marker-bits}. Hence $(f_h(v),f_h(w))\in E(H)$.
Thus $f_h:J_h(B)\to H$ is a homomorphism, proving necessity.

Conversely, assume the right side of~\eqref{neg:age-formula} and choose
$f_h:J_h(B)\to H$ for each $h\in B$. Retain the whole tables on $A_H$
and $B$. Set all old--new $E$-pairs true and every other mixed binary value
false. For $v\in D_h(B)$ prescribe
\begin{equation}\label{neg:completion-states}
 \begin{aligned}
 S_0(a_l,h,v)&\Longleftrightarrow l<f_h(v),\\
 S_1(a_l,h,v)&\Longleftrightarrow l=f_h(v),\\
 S_2(a_l,h,v)&\Longleftrightarrow l>f_h(v).
 \end{aligned}
\end{equation}
\begin{equation}\label{neg:completion-bits}
 \begin{aligned}
 U_L(b_{ij},h,v)&\Longleftrightarrow f_h(v)=i,\\
 U_R(b_{ij},h,v)&\Longleftrightarrow f_h(v)=j.
 \end{aligned}
\end{equation}
Every other mixed ternary value is false. These prescriptions concern
distinct ordered tuples. In particular, the choices for hubs $h$ and $v$
do not conflict. Every prescribed ternary tuple uses two distinct new
vertices, so each new point still has type $t_H$ over $A_H$.

It remains to check legality. The unary tables force every $O$-vertex of a
potential forbidden configuration into $A_H$ and every $V$-vertex into $B$.
The $E,Q$ hypotheses force each relevant $v$ to belong to $D_h(B)$.
Exactly one of $l<f_h(v)$, $l=f_h(v)$ and $l>f_h(v)$ holds, so
\eqref{neg:completion-states} assigns exactly one state. Since
$1\leq f_h(v)\leq k$, the endpoint states are $0$ and $2$, verifying
\eqref{neg:states}. Adjacent states belong to $\mathcal T$, which gives
\eqref{neg:transition}. The definitions also give
\[
 L(b_{ij},a)\Rightarrow a=a_i,\qquad
 R(b_{ij},a)\Rightarrow a=a_j;
\]
together with~\eqref{neg:completion-states}--\eqref{neg:completion-bits},
these implications give~\eqref{neg:copy-bits}.
Finally, a violation of~\eqref{neg:conflict} at $b_{ij},h,v,w$ would give
\[
 (v,w)\in E(J_h(B)),\qquad f_h(v)=i,\quad f_h(w)=j,
 \qquad(i,j)\in Z_H,
\]
contradicting that $f_h$ is a homomorphism. Thus the constructed $Y$ is
legal, and~\eqref{neg:completion} proves sufficiency.
\end{proof}

\begin{corollary}\label{neg:hom-order}
For finite nonempty loopless directed graphs $H,H'$,
\begin{equation}\label{neg:hom-order-equation}
 \mathcal C_H\subseteq\mathcal C_{H'}\quad\Longleftrightarrow\quad H\to H'.
\end{equation}
\end{corollary}

\begin{proof}
The reverse implication follows by composing the homomorphisms in
\eqref{neg:age-formula}. For the forward implication, form $B_H$ on
$\{h\}\dot\cup V(H)$ by giving all points type $\delta$, making $E$
complete, and setting
\[
 Q^{B_H}=\{(h,v):v\in H\},\qquad G^{B_H}=E(H).
\]
All other values are false. Then
$J_h(B_H)=H$ and $J_v(B_H)=\varnothing$ for $v\ne h$. Hence
$B_H\in\mathcal C_H$, and $B_H\in\mathcal C_{H'}$ implies $H\to H'$.
\end{proof}

\subsection{Infinite degrees}
\label{neg:infinite}

\begin{proof}[Proof of Theorems~\ref{thm:negative} and~\ref{neg:main}]
For $m\ge2$, let $D_m$ be the directed cycle on $\mathbb Z/m\mathbb Z$
with arcs $i\to i+1$. A homomorphism $g:D_p\to D_q$ necessarily satisfies
$g(i)=g(0)+i\pmod q$. Therefore
\[
 D_p\to D_q\quad\Longleftrightarrow\quad q\mid p.
\]
The ages $\mathcal C_{D_p}$, with $p$ ranging over the primes, are pairwise
incomparable by Corollary~\ref{neg:hom-order}. Their types $t_{D_p}$ all
have singleton type $\delta$ in the same structure $K_*$. For every
integer $q\geq2$, Proposition~\ref{neg:persistence} gives a $q$-colouring of
$(K_*)_\delta$ using all $q$ colours in every copy of $K_*$. Thus
\[
 \brd(\delta,K_*)=\crd(\delta,K_*)=\infty.
\]

Let $A\in\Age(K_*)$ have $k\ge1$ vertices of type $\delta$, and fix one
such vertex $a\in A$. In every copy $K'\cong K_*$, each vertex of type
$\delta$ is the image of $a$ under an embedding of $A$ into $K'$.
For the preceding $q$-colouring $\chi_q$, the colouring
\[
 e\longmapsto\chi_q(e(a)),\qquad e\in\Emb(A,K_*),
\]
uses all $q$ colours on every copy, proving $\brd(A,K_*)=\infty$.

For copies of $A$, use the finite colouring
\[
 B\longmapsto\{\chi_q(x):x\in B\cap(K_*)_\delta\}.
\]
Each value has size at most $k$, and the union of its values over any
copy of $K_*$ contains all $q$ colours. A bound of $t$ values would imply
$q\le kt$. Choosing $q=kt+1$ excludes every finite $t$.
\end{proof}

\begin{corollary}\label{neg:no-brs}
The structure $K_*$ has no precompact relational big Ramsey structure.
More strongly, it has no precompact relational expansion in which every
realized expanded singleton has big Ramsey degree one.
\end{corollary}

\begin{proof}
Suppose such an expansion has $m$ singleton types above $\delta$.
Take the persistent colouring $\chi_{m+1}$ above. Successively apply
singleton degree one to the pullback of this colouring on each expanded
singleton type. Composing the resulting self-embeddings gives an expanded
copy on which each of the $m$ types is monochromatic. Its reduct is a copy
of $K_*$ with at most $m$ colours on the $\delta$-vertices, a contradiction.
\end{proof}

\begin{remark}\label{neg:scope}
The counterexample uses several unary, binary and ternary relations; its
ternary relations need not be symmetric. It is not a counterexample in the language
of a single symmetric ternary relation. In particular, the proof does not
settle the question about the homogeneous $3$-uniform hypergraph omitting
the irreducible four-vertex hypergraph with three edges. Nor do we assert
that every irreducible forbidden family containing a four-vertex structure
gives infinite degrees.
\end{remark}

\section*{AI Assistance}

The research presented in this paper was carried out
with substantial assistance from OpenAI's GPT-6.1 Sol, GPT-6 Astra
and GPT-6 Pro. The author directed the
research process, evaluated the generated arguments,
and assumes full responsibility for the final results. Further details on the associated token consumption are available at \url{https://github.com/DJfrankkk/Big-Ramsey-Structures-WWWith-token-usage}.

\end{document}